\documentclass[11pt,letterpaper]{amsart}

\usepackage{tikz}
\tikzstyle{every node}=[circle, draw, fill=black!50,
                        inner sep=0pt, minimum width=3pt]

\usepackage{amsfonts,amsmath,latexsym,color,epsfig,hyperref,enumitem, amssymb,bbm}

\newtheorem{theorem}{Theorem}
\newtheorem{lemma}{Lemma}

\newtheorem{cor}[lemma]{Corollary}
\newtheorem{prop}[lemma]{Proposition}

\renewcommand{\le}{\leqslant}
\renewcommand{\ge}{\geqslant}
\renewcommand{\leq}{\leqslant}
\renewcommand{\geq}{\geqslant}

\def\qed{\ifvmode\mbox{ }\else\unskip\fi\hskip 1em plus 10fill$\Box$}

\input{epsf}

\makeatletter
\def\Ddots{\mathinner{\mkern1mu\raise\p@
\vbox{\kern7\p@\hbox{.}}\mkern2mu
\raise4\p@\hbox{.}\mkern2mu\raise7\p@\hbox{.}\mkern1mu}}
\makeatother

\def\R{\mathbb R}
\def\Z{\mathbb Z}
\def\Q{\mathbb Q}
\def\C{\mathbb C}
\def\F{\mathbb F}

\def\N{\mathbb N}
\def\E{\mathbb E}

\newtheorem{innercustomgeneric}{\customgenericname}
\providecommand{\customgenericname}{}
\newcommand{\newcustomtheorem}[2]{%
  \newenvironment{#1}[1]
  {%
   \renewcommand\customgenericname{#2}%
   \renewcommand\theinnercustomgeneric{##1}%
   \innercustomgeneric
  }
  {\endinnercustomgeneric}
}

\newcustomtheorem{customthm}{Conjecture}
\newcustomtheorem{customlemma}{Lemma}

\title{\vspace{-0.7cm}Generalized arithmetic Kakeya}
\author{Cosmin Pohoata \and Dmitrii Zakharov}
\date{}

\address{Department of Mathematics, Emory University, Atlanta, GA.}
\email{cosmin.pohoata@emory.edu} 

\address{Department of Mathematics, Massachusetts Institute of Technology, Cambridge, MA 02139, USA}
\email{zakhdm@mit.edu}

\begin{document}

\begin{abstract}
Around the early 2000-s, Bourgain, Katz and Tao introduced an arithmetic approach to study Kakeya-type problems. They showed that the Euclidean Kakeya conjecture follows from a natural problem in additive combinatorics, now referred to as the `Arithmetic Kakeya Conjecture'. 

We consider a higher dimensional variant of this problem and prove an upper bound using a certain iterative argument.
The main new ingredient in our proof is a general way to  strengthen the sum-difference inequalities of Katz and Tao which might be of independent interest. As a corollary, we obtain a new lower bound for the Minkowski dimension of $(n,d)$-Besicovitch sets.
\end{abstract}

\maketitle

\section{Introduction}

The Kakeya problem asks for the smallest possible dimension that a set of points $\mathbb{R}^{n}$ may have if it contains a unit segment in every direction. Sets with this property are called Besicovitch sets. There are actually multiple questions here, depending on which notion of dimension one considers. The so-called Kakeya conjecture asserts that such sets have full (Hausdorff) dimension for every $n$ and it is related to important questions in harmonic analysis, PDE’s, combinatorics, and number theory \cite{Tao1, Tao2002}. The conjecture is known in the plane \cite{Davies} but widely open for any $n \ge 3$. 
We refer to \cite{KatzZahl, WangZahl2022, WangZahl2024} for recent progress in $\R^3$, and to \cite{Hickman, Zahl} for a recent new approach in higher dimensions. See \cite{KatzTao} for a survey of older results.

The following `Arithmetic Kakeya Conjecture' conjecture, attributed to Katz and Tao \cite{Tao2002}, would imply the Kakeya conjecture for Minkowski dimension (see \cite{Bourgain2}). For $r \in \Q$ let $\pi_r: \R^2 \rightarrow \R$ denote the projection $\pi_r(x,y) = x+ry$ and let $\pi_\infty(x,y) = y$.

\begin{customthm}{A}\label{ArithmeticKakeya}
    For any $\varepsilon >0$ there exists a finite set $R \subset \Q$ such that the following holds. Let $G \subset \Z^2$ be a finite set such that the map $\pi_\infty$ is injective when restricted on $G$. Then we have 
    \[
    |G| \le \max_{r \in R} |\pi_r(G)|^{1+\varepsilon}.
    \]
\end{customthm}

It has proven to be very fruitful to study the entropic versions of problems in additive combinatorics \cite{GreenMannersTao, Tao4}, see e.g. the recent resolution of the Polynomial Freiman--Ruzsa conjecture over finite fields \cite{GGMT2023, GGMT2024}. In \cite{GR}, Green and Ruzsa introduced several equivalent formulations of Conjecture \ref{ArithmeticKakeya} and in particular they showed that it is equivalent to the following statement about Shannon entropies of integer random variables.

\begin{customthm}{B}\label{KakeyaEntropy}
For any $\varepsilon >0$ there exists a finite set $R \subset \Q$ such that we have: 
\begin{equation}\label{ak}
H(Y) \le (1+\varepsilon) \sup_{r\in R} H(X + r Y),  
\end{equation}
for any integer random variables $X,Y$.
\end{customthm}

Here and in what follows, all random variables are assumed to be discrete and finitely supported.
See Section \ref{section2} for relevant definitions. In the setting of arithmetic Kakeya problems, the set-theoretic and entropic languages are completely equivalent. So one can always use the language which is more convenient at any given moment.

To make progress on Conjecture \ref{KakeyaEntropy} and, thus, on the Euclidean Kakeya conjecture, one may try to make the ``$\varepsilon$'' in (\ref{ak}) as small as possible (while having the flexibility of choosing an arbitrarily large $R$). For a finite $R \subset \Q$ define $\beta(R)$ to be the minimum real number $\beta$ such that 
\begin{equation}\label{eq:entropy_kakeya}
H(Y) \le \beta \max_{r \in R} H(X+r Y)    
\end{equation}
holds for any integer random variables $X, Y$.
With this notation, Conjecture \ref{KakeyaEntropy} becomes equivalent to the assertion that the limit
\[
\beta_0 = \liminf_{R} \beta(R)
\]
taken over all finite $R \subset \Q$ is equal to $1$. We refer to \cite{KT1, KT2, Lemm} for bounds on $\beta(R)$ for specific choices of $R$. Katz and Tao \cite{KT2} showed that $\beta_0 \le \alpha$ where $\alpha = 1.675\ldots$ is the solution to the equation $\alpha^3-4\alpha+2=0$ and this bound has not been improved for over 20 years. 

In this paper, we are interested in a higher dimensional variant of the arithmetic Kakeya problem. Let $d\ge 1$ and let $R \subset \Q^d$ be a finite set. Then we define $\beta(R)$ as the smallest number such that for any integer random variables $X, Y_1, \ldots, Y_d$ we have
\begin{equation}\label{eq:GAPkakeya}
H(Y_1, \ldots, Y_d) \le \beta(R) \max_{(r_1,\ldots,r_d) \in R} H(X + r_1 Y_1 + \ldots + r_d Y_d).    
\end{equation}
Note that the quantity $\beta(R)$ is well-defined only if $R$ affinely spans $\Q^d$ and that $\beta(\phi(R)) = \beta(R)$ for any affine isomorphism $\phi:\Q^d\rightarrow \Q^d$. 
The Euclidean counterpart of the problem of estimating $\beta(R)$ for $R \subset \Q^d$ is the Kakeya problem for $d$-dimensional disks (instead of unit segments). Oberlin \cite{Oberlin2} observed that a certain `homogeneous' strengthening of arithmetic Kakeya estimates can be iterated to obtain bounds on Lebesgue measure of $(n,d)$-Besicovitch sets and the related maximal operators. 

The combinatorial version of this strengthening can be phrased as follows. For a finite set $R \subset \Q$ define $\beta_h(R)$ as the smallest number such that 
\begin{equation}\label{eq:entropy-kakeya-hom}
    H(Y) + (\beta_h(R)-1) H(X|Y) \le \beta_h(R) \max_{r\in R} H(X+r Y)
\end{equation}
for any integer random variables $X,Y$. So, comparing this to (\ref{eq:entropy_kakeya}), we add an extra term on the left hand side. Note that in the setting of Conjecture \ref{ArithmeticKakeya} one assumes that the mapping $(x,y) \mapsto y$ is injective on $G$, which, in the entropy language, means that the random variable $X$ is determined by the value of $Y$, i.e. $H(X|Y)=0$. So the extra term in (\ref{eq:entropy-kakeya-hom}) is not immediately useful for the original Kakeya problem. On the other hand, the fact that (\ref{eq:entropy-kakeya-hom}) is `homogeneous' (i.e. that the sum of coefficients on the left and right hand sides coincide) turns out to be important for the higher dimensional versions of the Kakeya problem. 
Oberlin used the `basic iteration' of Katz and Tao \cite{KT2} to (in particular) show that 
\begin{equation}\label{eq:homogen}
\liminf_{R \subset \Q} \beta_h(R) \le 1+\frac{1}{\sqrt{2}}
\end{equation}
holds. A version of this bound (suited for the Kakeya maximal operators) can then be iteratively applied to Kakeya-type problems for disks.

The first result of this paper states that the homogeneous (\ref{eq:entropy-kakeya-hom}) and the original (\ref{eq:entropy_kakeya}) versions of the arithmetic Kakeya problem are equivalent.

\begin{theorem}\label{thm:homogenious}
    Let $R \subset \Q$ be a finite set of size at least two. Then $\beta_h(R) = \beta(R)$.
\end{theorem}

This result is interesting because we do not know an exact value of $\beta(R)$ or $\beta_h(R)$ for any set $R \subset \Q$ of size at least three. The interesting direction of this theorem is showing that $\beta_h(R)\le \beta(R)$ holds, i.e. that the inequality (\ref{eq:entropy_kakeya}) implies an apparently stronger inequality (\ref{eq:entropy-kakeya-hom}) with an extra term on the left hand side. It might be interesting to investigate further which arithmetic entropy inequalities can be `upgraded' to their homogeneous form like this.
We should point out that in practice it seems to be the case that any proof of a bound of the form $\beta(R) \le \beta_0$ for some pair $(R, \beta_0)$ can be adapted with some effort to also prove that $\beta_h(R) \le \beta_0$ holds. On the other hand, the argument of Katz and Tao \cite{KT2} showing $\liminf_{R} \beta(R) \le \alpha = 1.675\ldots$ is quite involved and the assumption that $\pi_\infty$ is injective on $G$ is used heavily throughout. So it would require some considerable effort to rewrite that proof for the homogeneous version of the problem. Our Theorem \ref{thm:homogenious} provides such an extension for free. 

As a corollary, we prove the following bound for the generalized arithmetic Kakeya problem:
\begin{theorem}\label{prop:gapkakeya}
    Let $R \subset \Q$ be a finite set of size at least two and let $\beta=\beta(R)$. Then for any $d\ge 2$:
    \begin{equation*}
    \beta(R^d) \le d \frac{(\frac{\beta}{\beta-1})^d}{(\frac{\beta}{\beta-1})^d-1}.
    \end{equation*}
\end{theorem}

For example, using the simple fact $\beta(\{0,1\}) = 2$ we obtain
\[
\beta(\{0,1\}^d) \le d \frac{2^d}{2^d-1}.
\]
Compare this to the simple bounds which hold for any $R$ which affinely spans $\Q^d$:
\[
d \le \beta(R) \le d+1.
\]
Using the result of Katz and Tao \cite{KT2} we conclude that for any $\varepsilon >0$ there exists $R \subset \Q^d$ with 
\[
\beta(R) \le d \frac{(\frac{\alpha}{\alpha-1})^d}{(\frac{\alpha}{\alpha-1})^d-1}+\varepsilon,
\]
where $\alpha=1.675\ldots$ solves $\alpha^3-4\alpha+2=0$. It might be interesting to explore upper and lower bounds for the quantities $\beta(R)$ for various choices $R\subset \Q^d$ further.

It is perhaps also worth mentioning that in the finite fields setup such iterative arguments were also used recently by Dhar, Dvir, and Lund in \cite{DDL} to improve (and simplify) the results of Ellenberg and Erman \cite{EE} on the size of Furstenberg sets in $\mathbb{F}_{q}^{n}$.

One motivation to study the quantity $\beta(R)$ for $R\subset \Q^d$ with $d\ge 2$ is the Euclidean Kakeya problem for disks. For $n \ge d\ge 1$ say that $K \subset \R^n$ is an {\em $(n, d)$ Besicovitch set} if for any $d$-dimensional subspace $U \subset \R^n$ there is a translate $(x+U \cap B^n(0,1)) \subset K$. For instance, $(n, 1)$ Besicovitch sets are the subject of the Euclidean Kakeya conjecture. While Besicovitch \cite{Besicovitch, Besicovitch2} famously showed that there exist $(n,1)$-Besicovitch sets with zero Lebesgue measure, it is believed (\cite{Bourgain1991}, Remark 3.5) that already for $d= 2$ any $(n,d)$-Besicovitch set must have positive Lebesgue measure. 
In support of this belief, Bourgain \cite{Bourgain1991} showed that $(n,d)$ sets have positive measure as long as $d$ satisfies $2^{d}+d \geq n$. Some improvements were obtained by Oberlin \cite{Oberlin, OberlinThesis, Oberlin2}. In \cite{Oberlin2}, Oberlin used a version of the arithmetic approach of Katz and Tao to study $(n, d)$-Besicovitch sets and related maximal operators and showed that $(n,d)$-Besicovitch sets have positive measure as long as $(1+\sqrt{2})^{n-1}+n > d$.  

Theorem \ref{prop:gapkakeya} implies a bound on Minkowski dimension of $(n,d)$-Besicovitch sets.
\begin{cor}
    Let $K$ be a $(n,d)$-Besicovitch set. Then for any $R \subset \Q^d$:
    \[
    \dim_M K \ge \frac{d}{\beta(R)} n.
    \]
    In particular, we get 
    \[
    \dim_M K \ge \frac{(\frac{\alpha}{\alpha-1})^d-1}{(\frac{\alpha}{\alpha-1})^d} n,
    \]
    where $\alpha=1.675\ldots$ solves $\alpha^3-4\alpha+2=0$.
\end{cor}

In Section \ref{section2} we recall basic facts about Shannon entropy and some tools from additive combinatorics. 
In Section \ref{section3} we prove an auxiliary lemma. 
In Section \ref{section4} we prove Theorem \ref{thm:homogenious}. In Section \ref{section5} we prove Theorem \ref{prop:gapkakeya}. 

\subsection{Acknowledgements} We would like to thank Zeev Dvir for useful discussions. We also gratefully acknowledge the hospitality and support of the Institute for Advanced Study during the Special Year on Dynamics, Additive Number Theory and Algebraic Geometry, where a significant portion of this research was conducted. 

\section{Preliminaries}\label{section2}

\subsection{Shannon entropy.}
Let $X$ be a discrete random variable, the Shannon entropy $H(X)$ is defined by the formula
$$H(X) := \sum_{x} \Pr(X = x) \log \frac{1}{\Pr(X=x)}$$
where the summation is over all the elements $x$ in the range of $X$ for which $\Pr(X=x)$ is nonzero. Given two random variables $X$ and $Y$, define the conditional entropy as
\[
H(X|Y) := \sum_{y} \Pr(Y=y)H(X|Y=y).
\]
We have the following chain rule for the entropy of the joint random variable $(X,Y)$:
\[
H(X,Y) = H(X|Y) + H(Y) = H(Y|X) + H(X).
\]
More generally, for any choice of discrete random variables $X_{1},\ldots,X_{n}, Z$, we have the conditional chain rule:
\[
H(X_{1},\ldots,X_{n} | Z) = \sum_{i=1}^{n} H(X_i | X_1,\ldots,X_{i-1},Z).
\]
Another basic property is that conditioning on a random variable cannot increase entropy:
\[
H(X|Y) \leq H(X).
\]
An important consequence of this result is the submodularity of entropy:
\[
H(X,Y|Z) \leq H(X|Z) + H(Y|Z).
\]
Finally, we will need the Shearer’s
inequality. Here and throughout we use $[n]$ for $\left\{1,\ldots,n\right\}$.

\begin{lemma}[Shearer's inequality]\label{lem:shearer}
Let $\mathcal{F}$ be a family of subsets of $[n]$ (possibly with repetitions) with each $i \in [n]$ included in at least $t$ members of $\mathcal{F}$. For a random vector $(X_1,\ldots,X_n)$,
$$H(X_1,\ldots,X_n) \leq \frac{1}{t} \sum_{F \in \mathcal{F}} H(X_{F}),$$
where $X_{F}$ is the vector $(X_i: i \in F)$. 
\end{lemma}

\subsection{Ruzsa's equivalence and Freiman homomorphisms}

We will use a special case of Ruzsa's equivalence between sumsets and entropy \cite{Ruzsa}. Recall that $\pi_r: \Q^2 \rightarrow \Q$, for $r \in \Q\cup \{\infty\}$ is the linear projection given by $\pi_r(x, y) = x+ry$ and  $\pi_\infty(x, y) = y$.

\begin{lemma}\label{lem:ruzsa}
Let $X, Y$ be random variables on $\Z$ and $R \subset \Q$ be a finite set. Then for every $\varepsilon > 0$, there are infinitely many $n \ge 1$ for which there is a finite set $G \subset \Z^2$ such that:
\begin{itemize}
    \item[(i)] $|G| = 2^{(1\pm \varepsilon) H(X, Y) n}$,
    \item[(ii)] $|\pi_{r}(G)| = 2^{(1\pm \varepsilon) H(X+r Y)n}$ for $r \in R$,
    \item[(iii)] $|\pi_\infty(G)| = 2^{(1\pm \varepsilon) H(Y)n}$,
    \item[(iv)] $|G \cap \pi_\infty^{-1}(y)| = 2^{(1\pm \varepsilon) H(X|Y) n}$, for each $y \in \pi_\infty(G)$.
\end{itemize}
\end{lemma}

It will be convenient to take a Freiman embedding into a cyclic group.
The following is a variant of the standard Freiman embedding \cite[Lemma 5.26]{TaoVu}. Similar reduction to finite fields appears in the work of Green--Ruzsa \cite[Proof of Proposition 2.1]{GR}.

\begin{lemma}\label{lem:freiman}
    Let $G \subset \Z^2$ and $R \subset \Q$ be a non-empty finite set. Let $N= |\pi_\infty(G)|$ and suppose that $\pi_r(G) \le N$ for any $r \in R$ and $|\pi_\infty^{-1}(y) \cap G| \sim \frac{|G|}{N}$ for any $y \in \pi_\infty(G)$.
    Then for any prime $p \ge CN$ there is a subset $G' \subset (\Z/p\Z)^2$ such that
    $$
    |G'| \sim |G|,~ |\pi_\infty(G')| \sim |\pi_\infty(G)|,
    $$
    $$
    |\pi_r(G') | \lesssim |\pi_r(G)|, ~r\in R,
    $$
    where the implied constants depend on $R$ only.
\end{lemma}

\begin{proof}
Take a uniformly random element $\theta \in (0,1)$ and define a map $\phi_\theta: \Z \rightarrow \{0, \ldots, p-1\}$ by $\phi_\theta(x) = \lfloor p \{\theta x\}\rfloor$. We will take $G'$ to be a dense subset in $G'' = \{(\phi_\theta(x), \phi_\theta(y)),~ (x, y) \in G\}$ for an appropriately chosen $\theta$. 

Denote $B = \pi_\infty(G)$.
Say that $y \in B$ is bad if $\phi_\theta(y) = \phi_\theta(y')$ for some $y' \in B \setminus \{y\}$. For fixed pair $y\neq y'$ we have 
$$
\Pr[\phi_\theta(y) = \phi_\theta(y')] \le \Pr\left[ \|\theta (y-y')\|_{\R/\Z} \le \frac{1}{p} \right] \le \frac{2}{p}.
$$
So by union bound, $y$ is bad with probability at most $\frac{2|B|}{p} \le \frac{2}{C} \le 1/10$ if we take $C\ge 20$. For $y\in B$ let $G_y = \{x:~ (x, y) \in G\}$, by assumption, we have $|G_y| \sim \frac{|G|}{N}$ for any $y\in B$. Note that by assumption, for $r \in R$ we have $|\pi_r(G)| \le N$ holds. On the other hand, we have $|\pi_r(G)| \ge |G_y|$ for any $y \in B$. So each fiber has size at most $N$. Say that $(x, y) \in G$ is bad if there is $x'\neq x$ with $(x', y) \in G$ and $\pi_\theta(x') = \pi_\theta(x)$. 
By the union bound argument above, the probability that $(x, y)$ is bad is at most $\frac{1}{10}$. Define
$$
G_\theta = \{ (x,y)\in G:~ y \text{ is good and }(x,y)\text{ is good} \}.
$$
The expected size of $G_\theta$ is at least $0.8 |G|$, fix some $\theta \in (0,1)$ which realizes $|G_\theta| \ge 0.8|G|$ and define
$$
G' = \{(\phi_\theta(x), \phi_\theta(y)), ~ (x,y)\in G_\theta\} \subset \{0, \ldots, p-1\}^2.
$$
By design, we get $|G'| \ge 0.8 |G|$ and $|\pi_\infty(G')| \gtrsim N$ (by using that all fibers $G_y$ have approximately the same size). Let $r \in R$ and write $r = a/b$ for $a\in \Z$ and $b \in \N$. Then if $p$ is coprime to $a$ and $b$,
$$
|\pi_r(G')| = |\{ a \phi_\theta(x) + b \phi_\theta(y) \pmod p, ~ (x, y) \in G' \}|.
$$
For any $x, y\in \Z$ we have
$$
a \phi_\theta(x) + b \phi_\theta(y) - \phi_\theta(ax+b y) \in I 
$$
where $I =\{-|a|-b, \ldots, |a|+b\}\subset \Z/p\Z$. So we get
$$
\{ a \phi_\theta(x) + b \phi_\theta(y) \pmod p, ~ (x, y) \in G' \} \subset \{\phi_\theta(ax+by) \pmod p, ~ (x,y)\in G'\} + I
$$
which has cardinality upper bounded by $|\pi_r(G)| |I| \lesssim |\pi_r(G)|$. So the set $G'$ satisfies the claimed properties.
\end{proof}

\section{Sumsets with arithmetic progressions}\label{section3}

The following Lemma will be used in the proof of Theorem \ref{thm:homogenious}.

\begin{lemma}\label{lem_fou}
    Fix $\varepsilon > 0$ and let $p > p_0(\varepsilon)$ be a large enough prime and let $1\le m < p$. Let $A \subset \F_p$ be a set of size at least $p/m$ and let $I = \{1, \ldots, m\} \subset \F_p$. Then for at least $[m p^{-\varepsilon}]$ indices $j \in \{m/2, \ldots, m\}$ we have
    \begin{equation}\label{eq:union}
    |A + j^{-1} I| \ge p^{1-\varepsilon}.    
    \end{equation}
\end{lemma}

We recall some basic Fourier analysis over $\F_p$, which will be used in the proof of Lemma \ref{lem_fou}. For functions $f, g: \F_p \rightarrow \C$ we define the scalar product
$$
\langle f, g\rangle = \E_{x \in \F_p} f(x) g(x),
$$
convolution $f * g$:
$$
f * g(x) = \E_{y \in \F_p} f(y) g(x-y), 
$$
and the Fourier transform $\hat f: \F_p \rightarrow \C$: 
$$
\hat f(\xi) = \E_{x \in \F_p} e(-x\xi) f(x), ~~e(x) = e^{\frac{2\pi i x}{p}}.
$$

\begin{proof}
    By removing some elements from $A$ if necessary, we may assume that $|A| = \lceil \frac{p}{m} \rceil$.
    We need to show that $|A + j^{-1} I| \ge p^{1-\varepsilon}$ for many indices $j$. Let $J \subset \{m/2, \ldots, m\}$ be the set of primes between $m/2$ and $m$. For a pair of sets $A, B \subset \F_p$ we denote $E(A, B)$ the additive energy of the pair $A, B$, i.e. the number of quadruples $(a, a', b, b') \in A\times A\times B \times B$ such that $a+b=a'+b'$. Let $f = 1_A$ be the characteristic function of $A$.
    Let $d\ge C\varepsilon^{-1}$ be a large constant and let $\psi: I \rightarrow [0,1]$ be a bump function supported on $I$ such that 
    \begin{itemize}
        \item $\sum_{i \in I} \psi(i) \ge c_d |I|$, for some $c_d > 0$,
        \item $|\hat \psi(\xi)| \le C_d \frac{m}{p} (1+ \xi m/p)^{-d}$ for any $\xi \in [-p/2, p/2]$.
    \end{itemize}
    One can construct such $\psi$ by e.g. taking an iterated self-convolution of the characteristic function on a slightly smaller interval. Let $\psi_j(x) = \psi(j x)$, this is a bump function supported on the dilated interval $j^{-1} I$. 
    
    By the definition of the convolution, we have 
    $$
    \E_{x\in \F_p} f * \psi_j (x) = \E_{x, y\in \F_p} f(x) \psi_j(y) \ge c_d \frac{|A| |I|}{p^2}.
    $$
    On the other hand, by the Cauchy--Schwarz inequality, we have
    $$
    \E_{x\in \F_p} f * \psi_j (x) = \langle f * \psi_j, 1_{A + j^{-1} I} \rangle \le \|f * \psi_j \|_2 \|1_{A+j^{-1} I}\|_2.
    $$
    After some rearranging and noticing that $\|1_{S}\|_2 = |S|^{1/2} p^{-1/2}$ for any $S \subset \F_p$, this leads to a lower bound
    \begin{equation}\label{eq:lower}
        |A+ j^{-1} I| \ge c_d^2 \frac{|A|^2 |I|^2}{\|f* \psi_j\|_2^2 p^3}.
    \end{equation}
    Let us now study the following sum:
    $$
    E = \sum_{j \in J} \|f * \psi_j\|_2^2,
    $$
    which has the following Fourier expansion:
    $$
    E = \sum_{\xi \in \F_p} |\hat f(\xi)|^2 \sum_{j\in J} |\hat \psi_j(\xi)|^2. 
    $$
    Note that $\hat \psi_j(\xi) = \hat \psi(j^{-1}\xi)$. In particular, the $\xi = 0$ contribution is upper bounded by 
    $$
    |\hat f(0)|^2 |J| |\hat \psi(0)|^2 \le p^{-4} |A|^2 |J| m^2.
    $$
    Now fix $\xi \neq 0$ and let $K = p^{\varepsilon/2}$. Let $J_\xi$ be the set of $j \in J$ such that $j^{-1} \xi \in [-K \frac{p}{m}, K \frac{p}{m}]$. If we think of $\xi$ as an integer in the interval $[-p/2, p/2]$ then this condition holds if and only if $\xi = ij + k p$ for some integer $i \in [-K \frac{p}{m}, K \frac{p}{m}]$ and $k \in \Z$. Note that $j \in [1, m]$, so there are at most $2K$ options for $k$. For each fixed $k$, we have $\xi - k p = ij$ in $\Z$ and since $j$ is assumed to be prime between $m/2$ and $m$, there are at most $\frac{\log 2K p}{\log m/2} \le C\varepsilon^{-1}$ choices for $j \in J$ (coming from the prime factorization of $\xi-kp$). This gives an upper bound $|J_\xi| \le C\varepsilon^{-1} K$. Now we can estimate using the Fourier decay of $\psi$:
    $$
    \sum_{j \in J} |\hat \psi(\xi/j)|^2 \le C^2_d \left(|J_\xi| \frac{m^2}{p^2} + |J| K^{-2d} \frac{m^2}{p^2} \right) \le C_{\varepsilon, d} K \frac{m^2}{p^2}
    $$
    where we used $K^{-d} = p^{-\varepsilon d} < p^{-10}$, to get rid of the second term. Finally, using the identity $\sum_{\xi\in \F_p} |\hat f(\xi)|^2 = \frac{|A|}{p}$, we can bound $E$:
    $$
    E \le p^{-4} |A|^2 |J| m^2 + \sum_{\xi \in \F_p\setminus 0} |\hat f(\xi)|^2 C_{\varepsilon, d} K \frac{m^2}{p^2} \le p^{-4} |A|^2 |J| m^2 + C_{\varepsilon, d} K p^{-3} |A| m^2 \le 
    $$
    $$
    \le 2 p^{-2} |J| + 2C_{\varepsilon, d} K p^{-2} m \le 3C_{\varepsilon, d} K p^{-2} m
    $$
    where we used that $|A| = \lceil \frac{p}{m} \rceil$ and $|J| \le m$. By Markov, at least half of the indices $j \in J$ satisfy $\|f * \psi_j\|_2^2 \le \frac{2}{|J|} E$. So for these indices we obtain 
    $$
    |A + j^{-1} I| \ge c_d^2 \frac{|A|^2 m^2}{\frac{2}{|J|} E p^3} \ge c_{\varepsilon, d} K^{-1} \frac{|J|}{m} \frac{|A|^2 m^2}{p}
    $$
    We have $|A| m \ge p$, $K = p^{\varepsilon/2}$ and $|J| \ge c \frac{m}{\log m}$ by the Prime Number Theorem. So provided that $p$ is large enough, this implies $|A+j^{-1} I| \ge p^{1-\varepsilon}$ for at least half of the indices $j \in J \subset \{1, \ldots, m\}$. This completes the proof.
\end{proof}

\section{Proof of Theorem \ref{thm:homogenious}}\label{section4}

It is clear that $\beta(R) \le \beta_{h}(R)$ for any $R \subset \Q$, in the rest of the proof we show that $\beta(R) \ge \beta_h(R)$ holds as well.

We start with an arbitrary set $R \subset \Q$ and a pair of random variables $X,Y$. Fix an arbitrarily small $\varepsilon_0>0$. We are going to show the estimate 
\begin{equation}\label{eq:goal}
    (1-\varepsilon_0) H(Y) + (\beta(R)-1) H(X|Y) \le \beta(R) \max_{r \in R} H(X+rY).
\end{equation}
Taking $\varepsilon_0$ to approach zero would then complete the proof.

We may assume that $H(Y) \ge (1+\varepsilon_0) \max_{r \in R} H(X+rY)$ for any $r \in R$ since otherwise (\ref{eq:goal}) holds automatically. By replacing $X$ with $m X$ and $R$ with $m R$ for a constant $m = m(R)$ depending on $R$, we may assume that $R \subset \Z$ (let $m$ to be the product of all denominators of all $r \in R$). Thus, we may assume that $R \subset [-k, k] \cap \Z$ for some constant $k$ depending on $R$.

We perform a sequence of reductions. Apply Lemma \ref{lem:ruzsa} to the pair $(X, Y)$, the set $R$ and some parameter $\varepsilon \ll \varepsilon_0$. Let $G \subset \Z^2$ be the resulting graph. Next, apply Lemma \ref{lem:freiman} to the graph $G$ to obtain a new graph $G' \subset \F_p^2$ with the following properties:
\begin{itemize}
    \item $|G'| = 2^{(1\pm \varepsilon) H(X, Y) n}$,
    \item $|\pi_\infty(G')| \sim p = 2^{(1\pm \varepsilon) H(Y) n}$,
    \item $|\pi_r(G')| \le 2^{(1\pm \varepsilon) H(X+r Y)}$ for any $r \in R$,
    \item $|G' \cap \pi_\infty^{-1}(y)| = 2^{(1\pm \varepsilon) H(X|Y) n}$, for each $y \in \pi_\infty(G)$.
\end{itemize}

Let $\alpha = \varepsilon_1 + \max_{r \in R} \frac{H(X+rY)}{H(Y)}$ and $\gamma = \frac{H(X|Y)}{H(Y)} - \varepsilon_1$, where $\varepsilon_1=c\varepsilon_0$ for a sufficiently small constant $c>0$. Let $A = \bigcup_{r \in R} \pi_r(G')$ then we get an upper bound $|A| \le p^{\alpha-\varepsilon}$ provided that $\varepsilon$ is small enough compared to $\varepsilon_0$. 
Let $B = \pi_\infty(G')$, by the above, we have $|B| \sim p$. For $y \in B$ let $D_y = \{x:~ (x, y) \in G'\}$ be the fiber above $y$. We have for any $y\in B$
$$
|D_y| \ge 2^{(1\pm \varepsilon) H(X|Y) n} \ge p^{\gamma+\varepsilon}.
$$
Let $m = p^{1-\gamma}$, sample uniformly random elements $t \in \F_p$ and $d \in \F_p\setminus\{0\}$ and define random arithmetic progressions
\begin{align*}
I_0 &= \{d j, ~j = 0, \ldots, m-1\},\\
I &= I_0+t = \{t + d j, ~j = 0, \ldots, m-1\},\\
J &= \{d j, ~ j \in [m/2, m] \},
\end{align*}
and let us also put $\widetilde I = I+k J -kJ$ which is an arithmetic progression of length at most $2k m$ (recall that $k$ was chosen so that $R \subset [-k,k]$). Now we define a random set
$$
G_I' = \{(x, y) \in G':~ x \in I, ~ y \in J \}.
$$
Observe that if $(x, y) \in G_I'$ then for $r\in R$, 
$$
x+r y \in I + r J \subset \widetilde I,
$$
and so we get $\pi_r(G_I') \subset \widetilde I \cap A$. The size of the latter set can be easily controlled by Markov's inequality:
\begin{equation}\label{eq:markov}
\Pr\left[|\widetilde I \cap A| \ge T \frac{m}{p} |A|\right] \lesssim k T^{-1}.    
\end{equation}
Let us now understand the projection  $B_I =\pi_\infty(G_I')$. By definition $y \in B_I$ if and only if $y \in J$ and $D_y \cap I \neq \emptyset$. Note that $y\in J$ means that $d = y/j$ for some $j \in [m/2, m]$. Conditioned on $d = y/j$, we have $D_y \cap I \neq \emptyset $ if and only if
$$
t \in D_y - (y/j) I_0.
$$
So we conclude that 
\begin{equation}\label{eq:pry}
\Pr[y \in B_I] = \sum_{j \in [m/2, m]} \frac{|D_y- (y/j)I_0|}{p(p-1)}.    
\end{equation}
To lower bound this quantity, we apply Lemma \ref{lem_fou} to the set $A= -y^{-1} D_y$ and $\varepsilon>0$. By construction we have $|D_y| \ge p/m$ and so by Lemma \ref{lem_fou}, for at least $[m p^{-\varepsilon}]$ indices $j \in [m/2, m]$ we have $|D_y- (y/j)I_0| \ge p^{1-\varepsilon}$. So in (\ref{eq:pry}) we get
$$
\Pr[y \in B_I] \ge [m p^{-\varepsilon}] p^{-1-\varepsilon}
$$
Recall that $m = p^{1-\gamma}$, $\gamma = \frac{H(X|Y)}{H(Y)}-\varepsilon_1$. We have an estimate $H(X|Y) \le H(X+rY) \le H(Y)-\varepsilon_0$ and so it follows that $\gamma \le 1-c\varepsilon_0$. So we may assume that $m \ge p^{\varepsilon_1}$ and the lower bound on the probability above is non-trivial. We conclude that the expected size of $B_I$ is at least $m p^{-2\varepsilon}$. Note that $|B_I| \le m$ for any choice of $d, t$ (since $B_I \subset J$). So we have 
$$
\Pr[|B_I| \ge p^{-2\varepsilon}m] \ge p^{-2\varepsilon}.
$$
So if we take $T = p^{3\varepsilon}$ in (\ref{eq:markov}), then we obtain that with positive probability both inequalities $|B_I| \ge p^{-2\varepsilon}m$ and $|\widetilde I \cap A| \le p^{3\varepsilon-1}m |A|$ hold. Fix a pair $d, t$ which achieves this. We obtain a graph $G_I \subset I \times J$ such that
\begin{itemize}
    \item $\pi_\infty(G_I) \ge p^{-2\varepsilon} m$,
    \item $\pi_r(G_I) \le p^{3\varepsilon-1} m |A|$ for any $r \in R$.
\end{itemize}
Let $Y'$ be a uniformly random element in $\pi_\infty(G_I)$ and let $X'$ be a uniformly random element in the fiber $\pi_\infty^{-1}(Y') \cap G_I$. Then 
$$
H(Y') = \log |\pi_\infty(G_I)| \ge \log m - 2\varepsilon \log p =  (1-\gamma -2\varepsilon) \log p
$$ 
$$
H(X'+rY') \le \log |\pi_r(G_I)| \le (3\varepsilon-1)\log p +\log m + \log |A| \le (\alpha-\gamma +3\varepsilon) \log p
$$
and so by (\ref{eq:entropy_kakeya}) we obtain
$$
H(Y') \le \beta(R) \max_{r\in R} H(X'+rY'),
$$
$$
1-\gamma -2\varepsilon \le \beta(R) (\alpha-\gamma +3\varepsilon),
$$
recall that $\gamma = \frac{H(X|Y)}{H(Y)} -\varepsilon_1$ and $\alpha = \max_{r\in R} \frac{H(X+rY)}{H(Y)}+\varepsilon_1$. So rearranging leads to
$$
(1-\varepsilon_0) H(Y) - H(X|Y) \le \beta(R) (\max_{r\in R} H(X+r Y) - H(X|Y)).
$$
Taking $\varepsilon_0 \rightarrow 0$ then recovers the homogeneous Kakeya inequality (\ref{eq:entropy-kakeya-hom}). This completes the proof.

\section{Generalized Arithmetic Kakeya}\label{section5}

In this section we prove non-trivial bounds on $\beta(R)$ for $R \subset \Q^d$ in terms of the original arithmetic Kakeya inequalities. For $d\ge 1$ and $R \subset \Q^d$ let us denote $\delta(R) = \frac{\beta(R) / d}{\beta(R)/d-1}$, this quantity makes for cleaner statement of the result below.

\begin{prop}\label{thm:higher}
    Let $R_1, \ldots, R_d \subset \Q$ be finite sets and let $R = R_1 \times \ldots \times R_d$. Then we have
    \begin{equation*}
        \delta(R) \ge \delta(R_1) \delta(R_2) \ldots \delta(R_d).
    \end{equation*}
\end{prop}

Note that this implies Theorem \ref{prop:gapkakeya}.

\begin{proof}[Proof of Proposition \ref{thm:higher}]
    Let $X, Y_1, \ldots, Y_d$ be integer random variables. We want to prove an upper bound on $H(Y_1, \ldots, Y_d)$ in terms of $H(X+r_1Y_1+\ldots+r_dY_d)$ for $(r_1, \ldots,r_d) \in R = R_1\times \ldots \times R_d$. 
    We expand 
    \[
    H(Y_1, \ldots, Y_d) = \sum_{i=1}^d H(Y_i | Y_1, \ldots, Y_{i-1}).
    \]
    Fix an index $i \in \{1,\ldots, d\}$ and fix $r_{i+1} \in R_{i+1}, \ldots, r_d \in R_d$. We apply (\ref{eq:entropy-kakeya-hom}) to the variables 
    \begin{align*}
    \tilde X &= X+r_{i+1} Y_{i+1} + \ldots + r_d Y_d | Y_1, \ldots, Y_{i-1},\\
    \tilde Y &= Y_{i} | Y_1, \ldots, Y_{i-1}.
    \end{align*}
    and the set $R_i \subset \Q$. We obtain
    \[
    H(\tilde Y) + (\beta(R_i)-1) H(\tilde X | \tilde Y) \le \beta(R_i) \max_{r_i \in R_i} H(\tilde X + r_i \tilde Y) 
    \]
    which by expanding and denoting $\beta_i = \beta(R_i)$ gives 
    \begin{equation}\label{eq3}
    H(Y_i|Y_1, \ldots, Y_{i-1}) + (\beta_i-1) H(X+r_{i+1} Y_{i+1} + \ldots + r_d Y_d | Y_1, \ldots, Y_i) \le 
    \end{equation}
    \[
    \le \beta_i \max_{r_i \in R_i} H(X+r_{i} Y_{i} + \ldots + r_d Y_d | Y_1, \ldots, Y_{i-1}). 
    \]
    Let 
    \[
    H_i = \max_{r_{i+1} \in R_{i+1}, \ldots, r_d\in R_d} H(X+r_{i+1} Y_{i+1} + \ldots + r_d Y_d | Y_1, \ldots, Y_{i}),
    \]
    then (\ref{eq3}) gives
    \begin{equation}\label{eq4}
    H(Y_i|Y_1, \ldots, Y_{i-1}) \le \beta_i H_{i-1} - (\beta_i-1) H_i.    
    \end{equation}
    Let $\delta_i = \frac{\beta_i}{\beta_i-1}$, multiply both sides of (\ref{eq4}) by $\frac{1}{\delta_1 \ldots \delta_{i-1}\beta_i}$ 
    and sum over $i=1, \ldots, d$, we get:
    $$
    \sum_{i=1}^d \frac{1}{\delta_1 \ldots \delta_{i-1}\beta_i} H(Y_i|Y_1, \ldots, Y_{i-1}) \le  \sum_{i=1}^d \frac{\beta_i H_{i-1} - (\beta_i-1) H_i}{\delta_1 \ldots \delta_{i-1} \beta_i} = H_0 - \frac{H_d}{\delta_1\ldots \delta_d}.
    $$
    The sums of coefficients on both sides are the same since this is true for all instances of (\ref{eq:entropy-kakeya-hom}). If we expand $H(Y_i|Y_1, \ldots,Y_{i-1}) = H(Y_1, \ldots, Y_i) - H(Y_1, \ldots, Y_{i-1})$ and group same terms then one can check that the coefficient of $H(Y_1,\ldots, Y_i)$ in the resulting expression is equal to $\frac{\beta_{i+1}-\beta_i+1}{\delta_1\ldots\delta_{i-1}\beta_i\beta_{i+1}}$ which is a positive number since $\beta(R) \in (1,2]$ for any $R \subset \Q$. Now we can average this bound over all permutations of $Y_1, \ldots, Y_d$, apply Shearer's inequality (Lemma \ref{lem:shearer})
    \[
    H(Y_1, \ldots, Y_d) \le \frac{d}{k} \E_{i_1<\ldots< i_k} H(Y_{i_1}, \ldots, Y_{i_k})
    \]
    and conclude that
    \[
    \left(1-\frac{1}{\delta_1\ldots \delta_d}\right) \frac{1}{d} H(Y_1, \ldots, Y_d) \le H_0 - \frac{H_d}{\delta_1\ldots \delta_d} 
    \]
    \[
    \left(1-\frac{1}{\delta_1\ldots \delta_d}\right) \frac{1}{d} H(Y_1, \ldots, Y_d) + \frac{1}{\delta_1\ldots\delta_d} H(X|Y_1,\ldots Y_d) \le \max_{r_i\in R_i} H(X+r_1Y_1+\ldots r_dY_d)
    \]
    which after rearrangement gives us the desired bound. Note that we actually proved the stronger homogeneous version of (\ref{eq:GAPkakeya}).
\end{proof}

\end{document}